\documentclass{amsart}
\usepackage[all]{xy}
\usepackage{enumerate}
\usepackage{amssymb}
\usepackage{color}

\newtheorem{theorem}{\textbf{Theorem}}[section]
\newtheorem{lemma}[theorem]{\textbf{Lemma}}
\newtheorem{proposition}[theorem]{\textbf{Proposition}}
\newtheorem{proposition*}{Proposition}
\newtheorem{corollary}[theorem]{\textbf{Corollary}}

\theoremstyle{definition}
\newtheorem{definition}[theorem]{\textbf{Definition}}
\newtheorem*{example}{\textbf{Example}}

\newtheorem*{prf}{\textbf{Proof}}

\newtheorem*{remark}{\textbf{Remark}}
\newtheorem*{question}{\textbf{Question}}
\newtheorem*{remexa}{\textbf{Remarks and examples}}

\numberwithin{equation}{section}

\newcommand{\internalcomment}[1]{}

\newcommand{\N}{\mathbf{N}}
\newcommand{\Z}{\mathbf{Z}}

\newcommand{\ko}{\: , \;}

\newcommand{\ol}{\overline}

\renewcommand{\tilde}[1]{\widetilde{#1}}

\newcommand{\ra}{\rightarrow}

\newcommand{\arr}[1]{\stackrel{#1}{\longrightarrow}}

\newcommand{\opname}[1]{\operatorname{\mathsf{#1}}}

\newcommand{\Mod}{\opname{Mod}\nolimits}

\renewcommand{\ker}{\opname{ker}\nolimits}

\newcommand{\op}[1]{\opname{#1}\nolimits}

\newcommand{\cc}{{\mathcal C}}

\newcommand{\cf}{{\mathcal F}}

\newcommand{\ct}{{\mathcal T}}

\newcommand{\cx}{{\mathcal X}}
\newcommand{\cy}{{\mathcal Y}}

\renewcommand{\phi}{\varphi}

\newcommand{\Hom}{\opname{Hom}}

\newcommand{\Ext}{\opname{Ext}}

\begin{document}
\title{Classification of split torsion torsionfree triples in module categories}

\author{Pedro Nicol\'{a}s}

\address{Departamento de Matem\'{a}ticas Universidad de Murcia, Aptdo. 4021, 30100 Espinardo, Murcia, SPAIN}

\email{pedronz@um.es}

\author{Manuel Saor\'{\i}n}

\address{Departamento de Matem\'{a}ticas Universidad de Murcia, Aptdo. 4021, 30100 Espinardo, Murcia, SPAIN}

\email{msaorinc@um.es}
\thanks{The authors thank Ken Goodearl for telling them about key resuls in references \cite{Goodearl} and \cite{Small}. They
also thank the D.G.I. of the Spanish Ministry of Education and the Fundaci\'{o}n `S\'{e}neca' of Murcia for the financial support.}
\thanks{El segundo autor dedica este trabajo a su madre, Elisa Casta\~no, por su 75 cumplea\~nos.}


\begin{abstract}
A TTF-triple $(\cc,\ct,\cf)$ in an abelian category is
\emph{one-sided split} in case either $(\cc,\ct)$ or $(\ct,\cf)$
is a split torsion theory. In this paper we classify one-sided
split TTF-triples in module categories, thus completing Jans'
classification of two-sided split TTF-triples and answering a
question that has remained open for almost forty years.
\end{abstract}

\maketitle

\section{\textbf{Introduction}}
\bigskip

Since the early sixties torsion theories have played a important
role in Algebra. On one side, they Êtranslate to arbitrary module,
and in general abelian, Êcategories Êmany features of the
classical theory of torsion for abelian groups and modules over
PID's. On the other, they have been a fundamental tool to develop
a general theory of noncommutative localization, which imitates
the localization of commuative rings with respect to
multiplicative subsets. In modern times, the (pre)triangulated
version of them, namely, t-structures in triangulated categories
(cf. \cite{BeilinsonBernsteinDeligne} and \cite[Chapter
II]{BeligianisReiten}) are having a great impact in many fields
ranging from Representation Theory to Differential Geometry.

ÊIn the context of module categories, one of the concepts that has
deserved much attention is that of torsion torsionfree (TTF)
classes. It was introduced by Jans (\cite{Jans}), who gave a
bijection, for an arbitrary (associative unital) ring $A$, between
TTF-classes in the module category $\Mod A$ and idempotent
(two-sided) ideals of $A$ (cf. Proposition \ref{parametrizando
TTF-ternas} below). Idempotent ideals Ê(of not necessarily unital
rings) have recently had a great impact in  Homotopy Theory (cf.
\cite{Krause}) and  recollements of triangulated categories, which
are  the triangulated correspondent of TTF-triples, also play a
major role in different areas of Mathematics. That points in the
direction of a renewed interest in TTF-classes.

A TTF-class $\ct$ gives rise to a triple $(\cc,\ct,\cf )$, which
we shall call \emph{TTF-triple} in the sequel, where $(\cc ,\ct )$
and $(\ct ,\cf )$ are both torsion theories. Jans Êalso proved
that the above mentioned  bijection restricts to another one
between centrally split TTF triples in $\Mod A$ (see definition
below) and central idempotents of $A$. However, the existence of
TTF-triples for which only one of the torsion theories $(\cc ,\ct
)$ and $(\ct ,\cf )$ splits, which we shall call \emph{one-sided
split}, has been known for a long time (cf. \cite{Teply}) and no
classification of them has been available. That is, until now, the
idempotent ideals of $A$ which correspond by Jans' bijection to
those one-sided split TTF-triples have not been identified, even
though there were some efforts to classify those TTF-triples
Ê(see, e.g., \cite{Azumaya} and \cite{Ikeyama}). The goal of this
paper is to present such a classification, thus solving a problem
which has been open for almost forty years.

All rings appearing in the paper are associative with identity
and, unless explicitly said otherwise, all modules are right
modules, the category of which will be denoted $\Mod A$. Two-sided
ideals will be simply called ideals, Êwhen there is no risk of
confusion. The terminology that we use concerning Êrings and
modules is standard, and can be found in books like
\cite{AndersonFuller}, \cite{AuslanderReitenSmalo} or
\cite{Stenstrom}, and only in cases relevant to our work we shall
give precise definitions.

The organization of the paper is as follows. In section ~\ref{Torsion torsiofree triples} we give the definitions Êand relevant known results concerning
TTF-triples in a module category.
In section ~\ref{Left-split TTF-triples over arbitrary rings} we give the classification of left split TTF-triples in $\Mod A$
(cf. Theorem \ref{izquierda} and Corollary \ref{clasificacion-left-split}).
For the sake of clarity and as an intermediate step, we present in section ~\ref{Right-split TTF-triples over `good' rings} a partial classification of right split
TTF-triples, which is actually total
when the ring $A$ belongs to a class which includes semiperfect and left N\oe therian rings (cf. Theorem \ref{clasificacion-right-split-para-buenos} and Corollary
\ref{right-split para buenos2}). Then, Êin section ~\ref{Right-split TTF-triples over arbitrary rings}, Êwe classify the right split TTF-triples
$(\cc ,\ct ,\cf )$ in $\Mod A$ such that $A_A\in\cf$
(Theorem \ref{clasificacion-right-split}), from which the classification of all right split TTF-triples (Corollary \ref{clasification-right-split2}) follows.
\bigskip

\section{\textbf{Torsion torsiofree triples}}\label{Torsion torsiofree triples}
\bigskip

We refer the reader to Stenstr\"om's book \cite{Stenstrom} for the
terminology concerning torsion theories that we use in this paper.
We convene that if $(\cx ,\cy )$ is a torsion theory in $\Mod A$,
the associated idempotent radical, that we will call the
\emph{torsion radical}, will be denote by $\op{x}:\Mod
A\longrightarrow\Mod A$. Also, if  $\mathcal{Z}$ is a class of
$A$-modules, we shall put $\mathcal{Z}^{\bot}=\{X\in\Mod A\mid
Hom_A(Z,X)=0\text{ for all }Z\in\mathcal{Z}\}$ (resp.
$^{\bot}\mathcal{Z}=\{X\in\Mod A\mid Hom_A(X,Z)=0\text{ for all
}Z\in\mathcal{Z}\}$).

\begin{definition}
Let $A$ be an arbitrary ring. A \emph{torsion torsionfree triple} (or \emph{TTF-triple} for short) in $\Mod A$ is a triple $(\cc,\ct,\cf)$ formed by three
full subcategories of $\Mod A$ such that both $(\cc,\ct)$ and $(\ct,\cf)$ are torsion theories in $\Mod A$.
\end{definition}

The following result is due to Jans (\cite{Jans}, see also \cite[VI.8]{Stenstrom}).

\begin{proposition}\label{parametrizando TTF-ternas}
There exists a one-to-one correspondence between:
\begin{enumerate}[1)]
\item Idempotent Êideals $I^2=I$ of $A$.
\item TTF-classes in $\Mod A$ (i.e., full subcategories of $\Mod A$ closed under submodules, quotients, extensions and products).
\item TTF-triples in $\Mod A$.
\end{enumerate}
This correspondence maps the ideal $I$ to the full subcategory $\ct:=\{M\in\Mod A\mid MI=0\}$, the full subcategory $\ct$ to the
triple $(^{\bot}\ct,\ct,\ct^{\bot})$, and the TTF-triple $(\cc,\ct,\cf)$ to the ideal $\op{c}(A_{A})$.
\end{proposition}

We recall that if $(\cc ,\ct ,\cf )$ is the TTF-triple associated to the idempotent ideal $I$, then the torsion radicals associated
to $(\cc ,\ct )$ and $(\ct ,\cf )$ are given by $\op{c}(M)=MI$ and $\op{t}(M)=\op{ann}_M(I):=\{x\in M\mid xI=0\}$, respectively, for all $M\in \Mod A$.

Recall also that a torsion theory $(\cx,\cy)$ in $\Mod A$ splits if $\op{x}(M)$ is a direct summand of $M$, for every $A$-module $M$.

\begin{definition}
ÊLet $(\cc,\ct,\cf)$ be a TTF-triple in $\Mod A$. It will be called \emph{left-split} (resp. \emph{right-split}) if the torsion theory $(\cc,\ct)$ (resp. $(\ct,\cf)$)
Êsplits. It will be called \emph{centrally split} if it is both left and right split.
\end{definition}

The following is also well-known (cf. \cite[Proposition VI.8.5]{Stenstrom}).

\begin{proposition}
The one-to-one correspondence of Proposition ~\ref{parametrizando TTF-ternas} restricts to a one-to-one correspondence between:
\begin{enumerate}[1)]
\item Centrally split TTF-triples in $\Mod A$.
\item (Ideals of $A$ generated by) central idempotents of $A$.
\end{enumerate}
\end{proposition}

Put $\mathfrak{L}$, $\mathfrak{C}$ and $\mathfrak{R}$ for the sets of left, centrally and right split TTF-triples in $\Mod A$.

Since there are one-sided split TTF-triples which are not centrally split (cf. \cite{Teply}) Êwe should have a diagram of the form:
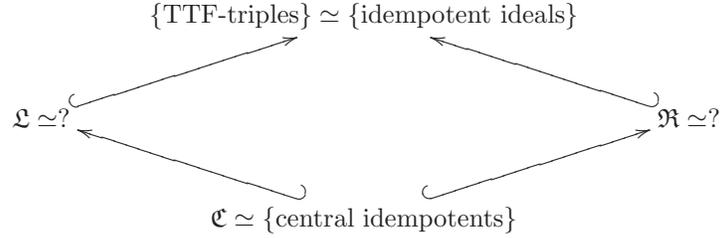
\[\xymatrix{ & \{\text{TTF-triples}\}\simeq\{\text{idempotent ideals}\} &  \\
\mathfrak{L}\simeq ?\ar@{^(->}[ur] & & \mathfrak{R}\simeq ?\ar@{_(->}[ul] \\
& \mathfrak{C}\simeq\{\text{central idempotents}\}\ar@{_(->}[ul]\ar@{^(->}[ur] &
}
\]

The following is the main question Êtackled in the paper:

\begin{question}
ÊWhat should replace the question marks in the diagram above?
\end{question}
\bigskip

\section{\textbf{Left-split TTF-triples over arbitrary rings}}\label{Left-split TTF-triples over arbitrary rings}
\bigskip

The following description of modules over triangular rings will be frequently used.

\begin{remark}
Let $B$ and $C$ be rings, Ê$M$ be a $B$-$C$-bimodule and
$A=\scriptsize{\left[\begin{array}{cc}C&0\\M&B\end{array}\right]}$
be the associated \emph{triangular matrix ring}. It is well- known
(cf. \cite[Chapter III]{AuslanderReitenSmalo}) that $\Mod A$ is
equivalent to a category, $\cc_{A}$, whose objects are triples
$(X,Y;\varphi)$ where $X\in\Mod C\ko Y\in\Mod B$ and
$\varphi\in\Hom_{C}(Y\otimes_{B}M,X)$. We shall often identify
$\Mod A$ with $\cc_{A}$. Notice that the assignments
$X\rightsquigarrow (X,0;0)$ and $Y\rightsquigarrow (0,Y;0)$ give
(fully faithful) embeddings $\Mod C\ra\cc_{A}\simeq\Mod A$ and
$\Mod B\ra\cc_{A}\simeq\Mod A$. We shall frequently identify $\Mod
C$ and $\Mod B$ with their images by these embeddings.
\end{remark}

We can already give the main result of this section. Recall that if $C$ is a ring and $M$ is a $C$-module, then $M$ is
\emph{hereditary} $\Sigma$-\emph{injective} in case every quotient of a direct sum of copies of $M$ is injective.

\begin{theorem}\label{izquierda}
Let $A$ be a ring and let $\ct$ be a full subcategory of $\Mod A$. The following assertions are equivalent:
\begin{enumerate}[1)]
\item $(^{\bot}\ct,\ct,\ct^{\bot})$ is a left-split TTF-triple in $\Mod A$.
\item There exists an idempotent $e\in A$ such that:
\begin{enumerate}[i)]
\item $(1-e)Ae=0$.
\item $NeA$ is a direct summand of $N$ for every $N\in\Mod A$.
\item $\ct=\{N\in\Mod A\mid N(1-e)=N\}$.
\end{enumerate}
\item There exists a ring isomorphism
\[A\arr{\sim}\scriptsize{\left[\begin{array}{cc}C&0\\M&B\end{array}\right]}
\]
where:
\begin{enumerate}[i)]
\item $M$ is a $B$-$C$-bimodule such that $M$ is hereditary $\Sigma$-injective in $\Mod C$.
\item $\ct$ identifies with $\{(X,0,0)\mid X\in\Mod C\}\simeq\Mod C$.
\end{enumerate}
\end{enumerate}
\end{theorem}
\begin{prf}
$(1\Rightarrow 2)$ Put $(\cc,\ct,\cf):=(^{\bot}\ct,\ct,\ct^{\bot})$. By hypothesis $I=\op{c}(A)=eA$, for some idempotent $e\in A$, Ê
and so $(1-e)A\cong \frac{A}{c(A)}$ belongs to $\ct$. Hence $(1-e)Ae=0$.
We know that $\op{c}(N)=NI=NeA$, which is then a direct summand of $N$, for every $N\in\Mod A$.
We also have $\ct=\{N\in\Mod A\mid Ne=0\}=\{N\in\Mod A\mid N(1-e)=N\}$.

$(2\Rightarrow 1)$ Of course, $\ct$ is a TTF-class (cf. \cite[VI.8]{Stenstrom}). Then it only remains to prove that the torsion theory $(^{\bot}\ct,\ct)$ splits.
First of all, since $eA\in^{\bot}\ct$ we get that $\op{Gen}(eA)\subseteq^{\bot}\ct$. On the other hand, decompose an arbitrary $N\in^{\bot}\ct$ as
$N=NeA\oplus N'$ where $N'e=0$. Hence $N'=N'(1-e)\in\ct$, and so $N'=0$ and Ê$N\in\op{Gen}(eA)$. Then $(^{\bot}\ct,\ct)=(\op{Gen}(eA),\ct)$, Ê
which is a split torsion theory by $(2.ii)$.

$(1=2\Rightarrow 3)$ Put $(\cc,\ct,\cf):=(^{\bot}\ct,\ct,\ct^{\bot})$. Take $C:=(1-e)A(1-e)\ko B:=eAe$ and $M:=eA(1-e)$. All the conditions of
$(3)$ are clearly satisfied except, perhaps, that $M_{C}$ is hereditary $\Sigma$-injective. Let us prove it. In case $eA(1-e)=0$ we are done, so assume that $eA(1-e)\neq 0$.
For an arbitrary set $\cx$ put
\[T:=eA(1-e)A^{(\cx)}\in\ct\ko D:=eA^{(\cx)}\in\cc\ko F:=\left(\frac{eA}{eA(1-e)A}\right)^{(\cx)}\in\cf.
\]
The short exact sequence
\[0\ra T\overset{i}{\hookrightarrow}D\ra F\ra 0
\]
is not split. Indeed, if it splits we would have $T\in\ct\cap\cc=\{0\}$, which contradicts the assumption $eA(1-e)\neq 0$. Take now a non-zero epimorphism
$p:T\twoheadrightarrow E$, and let us prove that $E$ is injective over $C$. By doing the pushout of $p$ and $i$ we get the commutative diagram
\[\xymatrix{0\ar[r] & T\ar[r]^{i}\ar@{->>}[d]_{p} & D\ar[r]\ar@{->>}[d] & F\ar[r]\ar@{=}[d] & 0 \\
0\ar[r] & E\ar[r]_{\mu} & V\ar[r]^{\pi} & F\ar[r] & 0 }
\]
Notice that $0\neq E\in\ct$ and $V\in\cc$ and, hence, the lower horizontal sequence does not split either. Suppose that $E$ is not injective over $C$. Then there
exists a non-split short exact sequence
\[0\ra E\ra T'\ra T''\ra 0
\]
in $\ct=\Mod C$. By doing the pushout of $E\ra T'$ and $E\arr{\mu}V$ we get the following commutative diagram with exact rows and columns
\[\xymatrix{ & 0\ar[d] & 0\ar[d] & & \\
0\ar[r] & E\ar[r]^{\mu}\ar[d] & V\ar[r]^{\pi}\ar[d] & F\ar[r]\ar@{=}[d] & 0 \\
0\ar[r] & T'\ar[r]\ar[d] & W\ar[r]\ar[d] & F\ar[r] & 0 \\
& T''\ar@{=}[r]\ar[d] & T''\ar[d] & & \\
& 0 & 0 & & }
\]
Since $V\in\cc$ and $T''\in\ct$ the central vertical short exact sequence splits. Therefore we can rewrite the diagram as follows
\[\xymatrix{ & 0\ar[d] & 0\ar[d] & & \\
0\ar[r] & E\ar[r]^{\mu}\ar[d] & V\ar[rr]^{\pi}\ar[d]_{\scriptsize{\left[\begin{array}{c}1\\ 0\end{array}\right]}} & & F\ar[r]\ar@{=}[d] & 0 \\
0\ar[r] & T'\ar[r]\ar[d] & V\oplus T''\ar[rr]^{\scriptsize{\left[\begin{array}{cc}\pi & 0\end{array}\right]}}\ar[d]^{\scriptsize{\left[\begin{array}{cc}0 & 1\end{array}\right]}} & & F\ar[r] & 0 \\
& T''\ar@{=}[r]\ar[d] & T''\ar[d] & & \\
& 0 & 0 & & }
\]
Now $T'\cong \ker\scriptsize{\left[\begin{array}{cc}\pi & 0\end{array}\right]}=\ker(\pi)\oplus T''\cong E\oplus T''$ and the
diagram is forced to be isomorphic to:

\[\xymatrix{ & 0\ar[d] & Ê& 0\ar[d] & & \\
0\ar[r] & E\ar[rr]^{\mu}\ar[d]_{\scriptsize{\left[\begin{array}{c}1\\0\end{array}\right]}} & & V\ar[rr]^{\pi}\ar[d]_{\scriptsize{\left[\begin{array}{c}1\\ 0\end{array}\right]}} & & F\ar[r]\ar@{=}[d] & 0 \\
0\ar[r] & E\oplus T''\ar[rr]_{\scriptsize{\left[\begin{array}{cc}\mu & 0 \\ 0 & 1\end{array}\right]}}\ar[d]_{\scriptsize{\left[\begin{array}{cc}0&1\end{array}\right]}} & & V\oplus T''\ar[rr]^{\scriptsize{\left[\begin{array}{cc}\pi & 0\end{array}\right]}}\ar[d]^{\scriptsize{\left[\begin{array}{cc}0 & 1\end{array}\right]}} & & F\ar[r] & 0 \\
& T''\ar@{=}[rr]\ar[d] & & ÊT''\ar[d] & & \\
& 0 & & 0 & & }
\]
and then the short exact sequence $0\ra E\ra T'\ra T''\ra 0$ splits, against the hypothesis.

$(3\Rightarrow 2)$ We identify $A=\scriptsize{\left[\begin{array}{cc}C & 0 \\ M & B\end{array}\right]}$. Taking
$e=\scriptsize{\left[\begin{array}{cc}0 & 0 \\ 0 & 1\end{array}\right]}$ we trivially have $(1-e)Ae=0$, and $\ct=\Mod C=\{N\in\Mod A\mid N(1-e)=N\}$.
It only remains to prove that $NeA$ is a direct summand of $N$ for each $N\in\Mod A$. Of course we have $N=NeA+N(1-e)A$.
By the proof of \cite[Theorem 2.5.]{AssemSaorin} we know that $NeA\cap N(1-e)A=NeA(1-e)A$, and so $NeA\cap N(1-e)A$ is generated by $eA(1-e)A=
\scriptsize{\left(\begin{array}{cc}0&0\\M&0\end{array}\right)}$, which is hereditary $\varSigma$-injective in $\Mod C$. Then
$NeA\cap N(1-e)A$ is injective in $\Mod C$, and so it induces a decomposition $N(1-e)A=NeA(1-e)A\oplus N'$ in $\Mod C\subseteq\Mod A$.
Therefore $N=NeA+N(1-e)A=NeA+NeA(1-e)A+N'=N'+NeA$. But $N'\cap NeA\subseteq N(1-e)A\cap NeA=NeA(1-e)A$, and so $N'\cap
NeA\subseteq N'\cap NeA(1-e)A=0$. Hence $N=N'\oplus NeA$.
\qed
\end{prf}

As a direct consequence of the theorem, the classification of left split TTF-triples in $\Mod A$ is at hand.

\begin{corollary} \label{clasificacion-left-split}
Let $A$ be a ring. The one-to-one correspondence of Proposition ~\ref{parametrizando TTF-ternas} restricts to a one-to-one correspondence between:
\begin{enumerate}[1)]
\item Left-split TTF-triples in $\Mod A$.
\item Two-sided ideals of $A$ of the form $I=eA$ where $e$ is an idempotent of $A$ such that $eA(1-e)$ is hereditary $\Sigma$-injective as a right $(1-e)A(1-e)$-module.
\end{enumerate}
\end{corollary}
\begin{prf}
If $e\in A$ is an idempotent such that $I=eA$ is a two-sided ideal, then $Ae\subseteq eA$ and, hence, $(1-e)Ae=0$. Now apply (the proof of) Theorem \ref{izquierda}.
\qed
\end{prf}
\bigskip

\section{\textbf{Right-split TTF-triples over `good' rings}}\label{Right-split TTF-triples over `good' rings}
\bigskip

In this section and the next hereditary perfect rings will play an important role. We gather some known properties of them which will be useful:

\begin{proposition}\label{ciclicos proyectivos implica proyectivo}
Let $R$ be a right perfect right hereditary ring. Then it is semiprimary, hereditary on both sides and the class of projective
$R$-modules (on either side) is closed under taking products.
\end{proposition}
\begin{prf}
From Ê\cite[Corollary 2 and Theorem 3]{Small} if follows that $R$ is semiprimary and hereditary (whence coherent) on both sides. Then apply \cite[Theorem 3.3]{Chase}.
\qed
\end{prf}

We start with some Êproperties of (right-split) TTF-triples which will be used in the sequel. The following one is due to Azumaya \cite[Theorem 6]{Azumaya}:

\begin{lemma}\label{hereditary igual a pure igual a condicion aritmetica}
Let $(\cc,\ct,\cf)$ be a TTF-triple in $\Mod A$, and $I=\op{c}(A)$ be its associated idempotent ideal. The torsion
theory $(\cc ,\ct )$ is hereditary if, and only if, $I$ is pure as a left ideal.
\end{lemma}

\begin{proposition}\label{primeras propiedades de los right-split}
Let $(\cc,\ct,\cf)$ be a right-split TTF-triple in $\Mod A$. Then $(\cc,\ct)$ is a hereditary torsion theory and $\cc\subseteq\cf$.
\end{proposition}
\begin{prf}
Applying \cite[Theorem VI.7.1.]{Stenstrom} to the split hereditary torsion theory $(\ct,\cf)$, we get that $\ct$ is closed under taking injective envelopes. But then
\cite[Proposition VI.3.2.]{Stenstrom} says that $(\cc,\ct)$ is hereditary, and consequently $\cc\subseteq\cf$ by \cite[Lemma VI.8.3]{Stenstrom}.
\qed
\end{prf}

We first want to classify the following type of right-split TTF-triples.

\begin{proposition} \label{anhadida}
Let $(\cc,\ct,\cf)$ be a right split TTF-triple in $\Mod A$ and let $I$ be its associated idempotent ideal. The following conditions are equivalent:
\begin{enumerate}[1)]
\item $\cc$ is closed for products.
\item $I=Ae$, for some idempotent $e\in A$.
\item Ê$A/I$ has a projective cover in $\Mod A$.
\item Ê$I$ is finitely generated as a left ideal.
\end{enumerate}
\end{proposition}
\begin{prf}
$(1\Leftrightarrow 2)$ ÊSince $(\cc ,\ct )$ is hereditary (cf. Proposition \ref{primeras propiedades de los right-split}), the class $\mathcal{C}$ is closed
for products if, and only if, $\cc$ is a $TTF$-class. Then Ê\cite[Theorem 3]{Azumaya} applies.

$(1\Leftrightarrow 3)$ ÊIt follows from Ê\cite[Theorem 8]{Azumaya}.

$(2\Leftrightarrow 4)$ Since $I$ is pure as a left ideal (cf. Lemma \ref{hereditary igual a pure igual a condicion aritmetica} and
Proposition \ref{primeras propiedades de los right-split}), $I$ is finitely generated on the left if, and only if, $_AI$ is a direct summand of $_AA$. Ê
\qed
\end{prf}

Recall that if $B$ is a ring, then a Ê$B$-module $P$ is called \emph{hereditary projective} (resp. \emph{hereditary} $\Pi$-\emph{projective}) in case every
submodule of $P$ (resp. of a direct product of copies of $P$) is projective. Recall also that a $B$-module $N$ is called \emph{FP-injective} in case $\Ext^1_B(?,N)$
vanishes on all finitely presented $B$-modules. The following type of modules will appear in our classification theorem:

\begin{proposition}\label{proposicion de hereditary pi-projective dual}
Let $B$ be a ring and $M$ be a left $B$-module. The following conditions are equivalent:
\begin{enumerate}[1)]
\item For every bimodule structure $_{B}M_{C}$ and every $X\in\Mod C$, the right $B$-module $\Hom_{C}(M,X)$ is hereditary projective.
\item There exists a bimodule structure $_{B}M_{C}$ such that, for every $X\in\Mod C$, the right $B$-module $\Hom_{C}(M,X)$ is hereditary projective.
\item If $S=\op{End}(_{B}M)^{op}$ and $Q$ is the minimal injective cogenerator of $\Mod S$, then $\Hom_{S}(M,Q)$ is hereditary $\Pi$-projective in $\Mod B$.
\item The character module $M^+=\Hom_{\Z}(M,\mathbf{Q}/\Z)$ is hereditary $\Pi$-projective in $\Mod B$.
\item $\op{ann}_B(M)=eB$ for some idempotent $e\in B$, $\ol{B}=B/\op{ann}_B(M)$ is a hereditary perfect ring and $M$ is FP-injective as a
left $\overline{B}$-module.
\end{enumerate}
When $B$ is an algebra over a commutative ring $k$, the above assertions are equivalent to:
\begin{enumerate}[6)]
\item If $Q$ is a minimal injective cogenerator of $\Mod k$, then $D(M):=\Hom_{k}(M,Q)$ is hereditary $\Pi$-projective in $\Mod B$.
\end{enumerate}
\end{proposition}
\begin{prf}
$(1\Rightarrow 2)$ Clear.

$(1\Rightarrow 3)$ Use the universal property of the product.

$(3\Rightarrow 4)$ If $E$ is an arbitrary injective cogenerator of $\Mod S$ then we have a section $E\ra Q^{\cx}$ for some set $\cx$.
Then $\Hom_{S}(M,E)$ is a direct summand of $\Hom_{S}(M,Q)^{\cx}$ and hence also $\Hom_{S}(M,E)$ is hereditary $\Pi$-projective. In
particular, this is true for $E=\Hom_{\Z}(S,\mathbf{Q}/\Z)$. By adjunction we get that $M^+$ is hereditary $\Pi$-projective.

$(4\Rightarrow 1)$ Let $C$ be a ring and $_{B}M_{C}$ a bimodule structure on $M$. Since $E=\Hom_{\Z}(C,\mathbf{Q}/\Z)$ is an injective cogenerator
of $\Mod C$, every right $C$-module embeds in a direct product of copies of $E$. So, it will be enough to prove that $\Hom_{C}(M,E)$ is hereditary
$\Pi$-projective, which is clear by adjunction.

$(2\Rightarrow 4)$ Take $X=\Hom_{\Z}(C,\mathbf{Q}/\Z)$ in $(2)$, and then use the universal property of the product and the hom-tensor adjunction.

$(5\Rightarrow 4)$ Since $\op{ann}_{B}(M)=eB$ with $e$ an idempotent of $B$, we get that Êa right $\ol{B}$-module $P$ is
projective over $B$ Êif, and only if, it is projective over $\ol{B}$. Notice that every submodule of a direct product of
copies of $M^+_B$ is annihilated by $\op{ann}_B(M)$. Then, there is no loss of generality in replacing $B$ by $\ol{B}$ and assuming that
$B$ is hereditary perfect and $_BM$ is FP-injective. The proof is in that case reduced to check that $M^+_B$ is projective. That
follows from \cite[Proposition 1.15]{Skljarenko} (see Proposition \ref{ciclicos proyectivos implica proyectivo}).

$(1=4\Rightarrow 5)$ Put $S=\op{End}(_{B}M)^{op}$ and let $\varphi: B\ra\op{End}(M_{S})$ the canonical morphism. By $(1)$,
$\op{End}(M_{S})$ is hereditary projective as a right $B$-module. Then the short exact sequence Êin $\Mod B$
\[0\ra\op{ann}_{B}(M)\hookrightarrow B\twoheadrightarrow\op{Im}(\varphi)\ra 0
\]
splits, and so $\op{ann}_{B}(M)=eB$, for some idempotent $e\in B$.

Now, as in the proof of $(5\Rightarrow 4)$, there is no loss of generality in replacing $B$ by $\ol{B}$ and assuming that $_BM$ is
faithful. Then $B_B$ is a submodule of $\op{End}(M_S)$, which is hereditary $\Pi$-projective by (1). Then every submodule
of a direct product of copies of $B_B$ is projective, which implies that $B$ is hereditary perfect by \cite[Theorem 3.3]{Chase}. But then,
by \cite[Proposition 1.15]{Skljarenko}, the fact that $M^+_B$ is projective Êimplies that $_BM$ is FP-injective.

Finally, the equivalence of $(6)$ with $(1)$-$(3)$ follows as the equivalence of $(4)$ with $(1)$-$(3)$ but replacing $\Z$ by $k$.
\qed
\end{prf}

\begin{definition}
A module $_{B}M$ satisfying the equivalent conditions of the last proposition will be said to \emph{have hereditary $\Pi$-projective dual}.
\end{definition}

We can now give the desired partial classification.

\begin{theorem} \label{clasificacion-right-split-para-buenos}
Let $A$ be a ring, $(\cc,\ct,\cf)$ be a TTF-triple in $\Mod A$ and $\op{c}(A)=I$ be the associated idempotent Êideal of $A$. The following assertions are equivalent:
\begin{enumerate}[1)]
\item $(\cc,\ct,\cf)$ is right split and $I$ is finitely generated as a left ideal.
\item There is an idempotent $e$ of $A$ such that $I=Ae$ and the left $(1-e)A(1-e)$-module $(1-e)Ae$ has a hereditary $\Pi$-projective dual.
\item There exists a ring isomorphism $A\cong\scriptsize{\left[\begin{array}{cc}C & 0 \\ M & B\end{array}\right]}$ such that $\ct$ gets identified with
$\{(0,Y;0)\mid Y\in\Mod B\}\simeq\Mod B$ and $_{B}M$ has a hereditary $\Pi$-projective dual.
\item There exist rings $B'\ko H\ko C$, where $H$ is hereditary perfect, Êand bimodules $_{H}M_{C}\ko _{B'}N_{H}$ such that:
\begin{enumerate}[i)]
\item $A\cong \scriptsize{\left[\begin{array}{ccc}C & 0 & 0\\ M & H & 0 \\ 0 & N & B'\end{array}\right]}$.
\item $\ct\cong\Mod\scriptsize{\left[\begin{array}{cc}H & 0 \\ N & B'\end{array}\right]}$.
\item $_{H}M$ is faithful and FP-injective.
\end{enumerate}
\end{enumerate}
\end{theorem}
\begin{prf}
$(2\Leftrightarrow 3)$ is left as an exercise.

\vspace*{0.3cm}
$(3\Rightarrow 4)$ Put $A=\scriptsize{\left[\begin{array}{cc}C & 0 \\ M & B\end{array}\right]}$.
By Proposition \ref{proposicion de hereditary pi-projective dual}(5), there exists an idempotent $e'$ of $B$ such that $\op{ann}_{B}(M)=e'B$, Ê
$\ol{B}:=B/\op{ann}_{B}(M)$ Êis a hereditary perfect ring and $_{\ol{B}}M$ is FP-injective.
We put $B':=e'Be'=Be'\ko H:=(1-e')B(1-e')\cong\ol{B}$ and $N:=e'B(1-e')$, and we get a ring isomorphism $B\cong
\scriptsize{\left[\begin{array}{cc}H & 0 \\ N & B'\end{array}\right]}$.
So, there is no loss of generality in assuming $B=\scriptsize{\left[\begin{array}{cc}H & 0
\\ N & B'\end{array}\right]}$ and $e'=\scriptsize{\left[\begin{array}{cc}0 & 0 \\ 0 &
1\end{array}\right]}$, so that $\op{ann}_{B}(M)=e'B=\scriptsize{\left[\begin{array}{cc}0
& 0 \\ N & B'\end{array}\right]}$. Representing now the left $B$-module $M$ as a triple
$(_{H}X,_{B'}Y;\varphi:N\otimes_{H}X\ra Y)$ we necessarily get $Y=0$. Then, abusing of
notation, Êwe can put $M=\scriptsize{\left[\begin{array}{cc}M
\\ 0\end{array}\right]}$. That gives the desired triangularization
\[\scriptsize{\left[\begin{array}{cc}C & 0 \\ M & B\end{array}\right]}\cong\scriptsize{\left[\begin{array}{ccc}C & 0 & 0\\ M & H & 0\\
0 & N & B'\end{array}\right]}.
\]

\vspace*{0.3cm}
Ê$(4\Rightarrow 3)$ Take $B=\scriptsize{\left[\begin{array}{cc}H & 0 \\ N
& B'\end{array}\right]}$ and $_{B}M_{C}=_{\scriptsize{\left[\begin{array}{cc}H & 0 \\ N &
B'\end{array}\right]}}\scriptsize{\left[\begin{array}{c}M \\
0\end{array}\right]}_{C}$.

\vspace*{0.3cm}

$(1\Rightarrow 2)$ From Proposition \ref{anhadida} we know that $I=Ae$, for some
idempotent $e\in A$. Then $eA(1-e)=0$ and we can identify $A$ with
$\scriptsize{\left[\begin{array}{cc}C & 0
\\ M & B\end{array}\right]}$, where $C=eAe\ko M=(1-e)Ae$ and $B=(1-e)A(1-e)$. Then
$I=\scriptsize{\left[\begin{array}{cc}C & 0
\\ M & 0\end{array}\right]}$ and $\ct=\{N\in\Mod A\mid NI=0\}=
\{N\in\Mod A\mid ÊN=N(1-e)\}$. This latter subcategory gets identified with $\{(0,Y;0)\mid Y\in\Mod B\}$.

We want to identify now $\op{t}(N)$ for every $N\in\Mod A$. If $N=(X_{C},Y_{B};\varphi:Y\otimes_{B}M\ra X)$ then $\op{t}(N)=
\op{ann}_{N}(I)=\{(0,y)\in N=X\oplus Y\mid \varphi(y\otimes
M)=0\}$. Now, using the usual adjunction we have an isomorphism
\[\Hom_{C}(Y\otimes_{B}M,X)\arr{\sim}\Hom_{B}(Y,\Hom_{C}(M,X))\ko \varphi\mapsto\varphi^t.
\]
Then $\varphi(y\otimes M)=0$ if and only if $\varphi^t(y)=0$. Then
$t(N)=(0,\ker(\varphi^t),0)$ and, Ê since it is a direct summand of $N=(X,Y;\varphi)$ by
hypothesis, we get that $\ker(\varphi^t)$ is a direct summand of $Y$ in $\Mod B$, for
every $\varphi\in\Hom_{C}(Y\otimes_{B}M,X)$. But then $\ker(\psi)$ is a direct summand of
$Y$ in $\Mod B$ for every $\psi\in\Hom_{B}(Y,\Hom_{C}(M,X))$. Since this is valid for
arbitrary $X_{C}\ko Y_{B}$, we can take Êfor $Y_{B}$ an arbitrary projective and we get
that every $B$-submodule of $\Hom_{C}(M,X)$ is projective. Then by Proposition
~\ref{proposicion de hereditary pi-projective dual} we get that $_{B}M$ has hereditary
$\Pi$-projective dual.

\vspace*{0.3cm}

$(3\Rightarrow 1)$ We identify $A$ with
$\scriptsize{\left[\begin{array}{cc}C & 0 \\ M &
B\end{array}\right]}$. Take
$e=\scriptsize{\left[\begin{array}{cc}1 & 0 \\ 0 & 0\end{array}\right]}$ and $I=Ae=\scriptsize{\left[\begin{array}{cc}C & 0 \\
M & 0\end{array}\right]}$. As we saw before, if $N=(X,Y;\varphi)\in\Mod A$, we have the identity $\op{t}(N)=(0,\ker(\varphi^t);0)$. But $\op{Im}(\varphi^t)$ is a
$B$-submodule of $\Hom_{C}(M,X)$, whence projective in $\Mod B$.
Thus, $\ker(\varphi^t)$ is a direct summand of $Y$ in $\Mod B$,
and so $\op{t}(N)$ is a direct summand of $N$ in $\Mod A$.
\qed
\end{prf}

\begin{corollary} \label{right-split para buenos2}
Let $A$ be a ring. Then the one-to-one correspondence of Proposition ~\ref{parametrizando
TTF-ternas} restricts to a one-to-one correspondence between:
\begin{enumerate}[1)]
\item Right-split TTF-triples in $\Mod A$ whose associated idempotent ideal $I$ is finitely generated on the left.
\item Two-sided ideals of the form $I=Ae$, where $e$ is an idempotent of $A$ such that $_{(1-e)A(1-e)}(1-e)Ae$ has a hereditary $\Pi$-projective dual.
\end{enumerate}

In particular, when $A$ satisfies either one of the two following conditions, the class $(1)$ above covers all the right split TTF-triples in $\Mod A$:

\begin{enumerate}[i)]
\item $A$ is semiperfect.
\item Every idempotent Êideal of $A$ which is pure on the left is also finitely generated on the left (e.g. $A$ left N\oe therian).
\end{enumerate}
\end{corollary}
\begin{prf}
Using Proposition \ref{anhadida}, under conditions $i)$ or $ii)$ the idempotent Êideal associated to a right split TTF-triple in $\Mod A$ is always finitely
generated on the left. With that in mind, the result is a direct consequence of the foregoing theorem.
\qed
\end{prf}
\bigskip

\section{\textbf{Right-split TTF-triples over arbitrary rings}}\label{Right-split TTF-triples over arbitrary rings}
\bigskip

Let $(\cc,\ct,\cf)$ be a ÊTTF-triple in $\Mod A$ Êwith associated idempotent ideal $I$. If $\op{lann}_A(I)=t(A_A)=(1-\varepsilon )A$, for some idempotent
$\varepsilon \in A$, then $\varepsilon A(1-\varepsilon )=0$ and $I\subseteq \varepsilon A\varepsilon$.

\begin{proposition}
The one-to-one correspondence of Proposition ~\ref{parametrizando
TTF-ternas} restricts to a one-to-one correspondence between:
\begin{enumerate}[1)]
\item Right-split TTF-triples in $\Mod A$.
\item Idempotent ideals $I$ of $A$ such that, for some idempotent $\varepsilon\in A$, one has that $\op{lann}_A(I)=(1-\varepsilon )A$ and the
TTF-triple in $\Mod\varepsilon A\varepsilon$ associated to $I$ is
right-split.
\end{enumerate}
\end{proposition}
\begin{prf}
According to our previous comments, Êour goal reduces to prove that if $(\cc ,\ct,\cf )$ is a TTF-triple such that $t(A_A)=(1-\varepsilon )A$ is a direct
summand of $A_A$, then it is right split if, and only if, the TTF-triple defined by $I$ in $\Mod\varepsilon A\varepsilon$ is also right split.

`Only if' part: Every right $\varepsilon A\varepsilon $-module $X$ can be viewed as a right $A$-module (by defining $X\cdot (1-\varepsilon )A=0$), and
then $X=\op{t}(X)\oplus F=\op{ann}_{X}(I)\oplus F$.

`If' part: Conversely, suppose that the TTF-triple in $\Mod\varepsilon A\varepsilon$ associated to $I$ is right split, and
put $C=\varepsilon A\varepsilon\ko B=(1-\varepsilon )A(1-\varepsilon )$ and $M=(1-\varepsilon)A\varepsilon$. As usual, we can identify $A$ with
$\scriptsize{\left[\begin{array}{cc}C & 0 \\ M & B\end{array}\right]}$, and in this case also $I$ with $\scriptsize{\left[\begin{array}{cc}I & 0 \\0 &
0\end{array}\right]}$, where $I$ is an idempotent Êideal of $C$ such that the TTF-triple in $\Mod C$ associated to $I$ is right
split. ÊNotice that, since $\op{lann}_A(I)=(1-\varepsilon )A$, we have $MI=0$.Ê Now, let $N=(X_{C},Y_{B};\varphi)$
be a right $A$-module. Then $\op{t}(N)=\op{ann}_{N}(I)=(\op{ann}_{X}(I),Y;\tilde{\varphi})$, where $\tilde{\varphi}$ is given by the decomposition of
$\varphi:Y\otimes_BM\stackrel{\tilde{\varphi}}{\longrightarrow}\op{ann}_X(I)\stackrel{j}{\hookrightarrow}
X$. Consider now a retraction $p:X\ra\op{ann}_{X}(I)$ in $\Mod C$ for the canonical inclusion $j$, which exists because the
TTF-triple induced by $I$ in $\Mod C$ is right split. Then $p\circ\varphi =p\circ j\circ\tilde{\varphi}=1\circ\tilde{\varphi}=\tilde{\varphi}$,
which implies that $\scriptsize{\left[\begin{array}{cc}p & 1\end{array}\right]}:(X,Y;\varphi)=N\ra t(N)=(\op{ann}_{X}(I),Y;\tilde{\varphi})$
is a morphism of $A$-modules, which is then a retraction for the canonical inclusion $t(N)\hookrightarrow N$.
\qed
\end{prf}

In the situation of the above proposition, one has that $\op{lann}_{\varepsilon A\varepsilon}(I)=0$, that is, the TTF-triple $(\cc',\ct ',\cf ')$ in
$\Mod\varepsilon A\varepsilon$ associated to $I$ has
the property that $\varepsilon A\varepsilon_{\varepsilon A\varepsilon}\in\cf'$. The problem of classifying right split TTF-triples gets then
reduced to answer the following:

\begin{question} Let $I$ be an idempotent Êideal of a ring $A$ such that $\op{lann}_{A}(I)=0$ (i.e., $A_{A}\in\cf$ where $(\cc,\ct,\cf)$ is the associated
TTF-triple in $\Mod A$). Which conditions on $I$ are equivalent to say that $(\cc,\ct,\cf)$ is right split?
\end{question}

Given a Êright $A$-module $M$ and a submodule $N$ , we shall say that $N$ is \emph{I-saturated} in $M$ when $xI\subseteq N$, with $x\in M$, implies that $x\in N$.
Equivalently, when $\frac{M}{N}\in\cf$.

When $X$ is a subset of $A$, we shall denote by $\mathcal{M}_{n\times n}(X)$ the subset of matrices of $\mathcal{M}_{n\times n}(A)$ with entries
in $X$.

\begin{lemma} \label{new-splitting}
Let $I$ be an idempotent ideal of the ring $A$. The following assertions are equivalent:
\begin{enumerate}[1)]
\item For every integer $n>0$ and every $\mathcal{M}_{n\times n}(I)$-saturated right ideal $\mathfrak{a}$ of $\mathcal{M}_{n\times n}(A)$,
there exists Ê$x\in\mathfrak{a}$ such that $(1_n-x)\mathfrak{a}\subseteq\mathcal{M}_{n\times n}(I)$.
\item For every integer $n>0$ and every $I$-saturated submodule $K$ of $A^{(n)}$, the quotient $\frac{A^{(n)}}{K+I^{(n)}}$ is
projective as a right $\frac{A}{I}$-module.
\end{enumerate}
\end{lemma}
\begin{prf}
Fix any integer $n>0$ and consider the equivalence of categories $F=:\Hom_A(A^{(n)},?):\Mod A\arr{\sim}\Mod \mathcal{M}_{n\times n}(A)$. It establishes a
bijection between $I$-saturated submodules $K$ of $A^{(n)}$ and $\mathcal{M}_{n\times n}(I)$-saturated right ideals $\mathfrak{a}$ of the ring
$\mathcal{M}_{n\times
n}(A)=\Hom_A(A^{(n)},A^{(n)})$. Now if $K$ and $\mathfrak{a}$ correspond by that bijection, one has that
$F(\frac{A^{(n)}}{K+I^{(n)}})\cong\frac{\mathcal{M}_{n\times n}(A)}{\mathfrak{a}+\mathcal{M}_{n\times n}(I)}$. Assertion $(2)$ is equivalent to say that,
for such a $K$, Êthe canonical projection $(\frac{A}{I})^{(n)}\cong\frac{A^{(n)}}{I^{(n)}}\twoheadrightarrow
\frac{A^{(n)}}{K+I^{(n)}}$ is a retraction Êin $\Mod A$. The proof is whence reduced to check that condition $(1)$ is equivalent to say that the canonical projection
$\frac{\mathcal{M}_{n\times n}(A)}{\mathcal{M}_{n\times
n}(I)}\twoheadrightarrow\frac{\mathcal{M}_{n\times
n}(A)}{\mathfrak{a}+\mathcal{M}_{n\times n}(I)}$ is a retraction in $\Mod \mathcal{M}_{n\times n}(A)$, for every
$\mathcal{M}_{n\times n}(I)$-saturated right ideal
$\mathfrak{a}$ of $\mathcal{M}_{n\times n}(A)$. To do that it is not restrictive to
assume that $n=1$, something that we do from now on in this proof. Then the existence of
an element $x\in\mathfrak{a}$ such that $(1-x)\mathfrak{a}\subseteq I$ is equivalent to
say that there is an element $x\in\mathfrak{a}$ such that $\bar{x}=x+I$ generates
$\frac{\mathfrak{a}+I}{I}$ and $\bar{x}^2=\bar{x}$. That is clearly equivalent to say
that the canonical projection $\frac{A}{I}\twoheadrightarrow\frac{A}{\mathfrak{a}+I}$ is
a retraction.
\qed
\end{prf}

\begin{definition}
An idempotent ideal $I$ of a ring $A$ will be called \emph{right
splitting} if it satisfies one (and hence both) of the equivalent
conditions of the above proposition, is pure as a left ideal and
$\op{lann}_{A}(I)=0$.
\end{definition}

\begin{example} \label{ejemplo de ideal right-splitting}
Let $I$ be an idempotent ideal of a ring $A$. If $I$ is pure on the left, $\op{lann}_{A}(I)=0$ and $A/I$ is a semisimple ring, then $I$ is right splitting.
\end{example}

\begin{theorem} \label{clasificacion-right-split}
Let $A$ be a ring and $(\cc,\ct,\cf)$ a TTF-triple in $\Mod A$ with associated idempotent ideal $I$ and such that $A_{A}\in\cf$.
The following assertions are equivalent:
\begin{enumerate}[1)]
\item $(\cc,\ct,\cf)$ is right split.
\item $(\cc,\ct)$ is hereditary and $F/FI$ is a projective right $A/I$-module for all $F\in\cf$.
\item $I$ is right splitting and $A/I$ is a hereditary perfect ring.
\end{enumerate}
\end{theorem}
\begin{prf}
$(1\Rightarrow 2)$ By Proposition ~\ref{primeras propiedades de los right-split} we know that $(\cc,\ct)$ is hereditary. Now take
arbitrary modules $T\in\ct\ko F\in\cf$, and apply $\Hom_{A}(?,T)$ to the short exact sequence
\[0\ra FI\hookrightarrow F\ra F/FI\ra 0.
\]
One gets the exact sequence
\[0=\Hom_{A}(FI,T)\ra\Ext^1_{A}(F/FI,T)\ra\Ext^1_{A}(F,T).
\]
The split condition of $(\ct,\cf)$ gives that $\Ext^1_{A}(F,T)=0$, and hence $\Ext^1_{A}(F/FI,T)=0$ for all $T\in\ct=\Mod\frac{A}{I}$.
Then $F/FI$ is a projective right $A/I$-module.

$(2\Rightarrow 1)$ Since $(\cc,\ct)$ is hereditary, it follows that $\op{t}(M)\cap MI=0$ for every $M\in\Mod A$. So the composition
\[\op{t}(M)\overset{j}{\hookrightarrow}M\overset{p}{\twoheadrightarrow}M/MI
\]
is a monomorphism with cokernel $\op{Cok}(p\circ j)\cong\frac{M}{\op{t}(M)+MI}\cong\frac{M/\op{t}(M)}{(M/\op{t}(M))\cdot I}$,
which, by hypothesis, is a projective right $A/I$-module. Then $p\circ j$ is a section in ($\Mod\frac{A}{I}$, and so in) $\Mod A$. Thus $j$ is a section in $\Mod A$.

$(2\Rightarrow 3)$ By Lemma ~\ref{hereditary igual a pure igual a condicion aritmetica} we know that $I$ is pure on the left.
Moreover, since $A_{A}\in\cf$ we get $\op{lann}_{A}(I)=0$. On the other hand, the fact that $F/FI$ is projective over $A/I$, for all
$F\in\cf$, implies that if $K<a^{(n)}$ is an $I$-saturated submodule then $\frac{A^{(n)}}{K+I^{(n)}}$ is projective as a
right $\frac{A}{I}$-module. Then, from Lemma \ref{new-splitting}, we derive that $I$ is right splitting.

Take now any right ideal $\mathfrak{b}/I$ of $A/I$ (notice that we then have $\mathfrak{b}\in\cf$). Since $(\cc,\ct)$ is hereditary,
then $\op{c}$ is left exact (cf. \cite[Proposition VI.3.1]{Stenstrom}, and so $\op{c}(\mathfrak{b})=\mathfrak{b}\cap\op{c}(A)$,
i.e. $\mathfrak{b}I=\mathfrak{b}\cap I=I$. Thus $\frac{\mathfrak{b}}{I}=\frac{\mathfrak{b}}{\mathfrak{b}I}$ is a projective right $A/I$-module.
That proves that $A/I$ is a right hereditary ring.

Finally, observe that if $\cx$ is any set then $A^{\cx}/A^{\cx}I$ is a projective right $A/I$-module. Now, \cite[Theorem 5.1]{Goodearl} says that $A/I$ is right perfect.

$(3\Rightarrow 2)$ By Lemma ~\ref{hereditary igual a pure igual a condicion aritmetica} we know that $(\cc,\ct)$ is hereditary. Let
us take now $F\in\cf$. If $F$ is finitely generated then $F\cong\frac{A^{(n)}}{K}$, for some $I$-saturated submodule
$K<a^{(n)}$. Then, by Lemma \ref{new-splitting}, we have that $F/FI\cong\frac{A^{(n)}}{K+I^{(n)}}$ is projective as a right
$A/I$-module. In case $F$ is not necessarily finitely generated, then Ê$F=\bigcup F_\alpha$ is the directed union of its finitely
generated submodules, which implies that $F/FI$ is a direct limit of the $F_\alpha/F_\alpha I$. Then $F/FI$ is a direct limit of
projective right $A/I$-modules. Since $A/I$ is right perfect, we conclude that $F/FI$ is projective over $A/I$, for all $F\in\cf$.
\qed
\end{prf}

The desired full classification of right split TTF-triples in $\Mod A$ is now available:

\begin{corollary} \label{clasification-right-split2}
Let $A$ be an arbitrary ring. The one-to-one correspondence of Proposition ~\ref{parametrizando TTF-ternas} restricts to a one-to-one correspondence between:
\begin{enumerate}[1)]
\item Right-split TTF-triples in $\Mod A$.
\item Idempotent ideals $I$ such that, for some idempotent $\varepsilon\in A$, $\op{lann}_A(I)=(1-\varepsilon )A$ and $I$ is a right splitting ideal of
$\varepsilon A\varepsilon $ with $\varepsilon A\varepsilon /I$ a hereditary perfect ring.
\end{enumerate}
\end{corollary}

\begin{remexa}
\begin{enumerate}[1)]
\item If in the definition of right splitting ideal we replace the conditions Êof Lemma \ref{new-splitting} by its corresponding ones with $n=1$, then
the correspondent of Theorem \ref{clasificacion-right-split} is not necessarily true. To see that, consider $A= \scriptsize{\left[\begin{array}{ccc}k & 0 \\ M & H
\end{array}\right]}$, where $k$ is an algebraically closed field, $H=k(\xymatrix{\cdot\ar@<0.75ex>[r]\ar@<-0.75ex>[r] & \cdot})$ is the Kronecker algebra
and $_HM=\tau_H (S)$, where $\tau_H$ is the Auslander-Reiten translation \cite[Chapter VII]{AuslanderReitenSmalo} and $S$ is the
simple injective left $H$-module. Notice that $_{H}M$ is faithful. We put $e=\scriptsize{\left[\begin{array}{ccc}1 & 0 \\ 0 & 0 \end{array}\right]}$ and take
$I=Ae=\scriptsize{\left[\begin{array}{ccc}k & 0 \\ M & 0 \end{array}\right]}$, which is clearly pure as a left ideal and satisfies that
$\op{lann}_A(I)=0$. A right ideal $\mathfrak{a}$ of $A$ is represented by a triple $(V,\mathfrak{u};\varphi:\mathfrak{u}\otimes_HM\longrightarrow V)$,
where $V$ is a vector subspace of $k\oplus M$, $\mathfrak{u}$ is a right ideal of $H$ such that $0\oplus\mathfrak{u}M\subseteq V$ and $\varphi$ is the
canonical multiplication map $h\otimes m\rightsquigarrow (0,hm)$. One readily sees that $\mathfrak{a}$ is $I$-saturated if, and only
if, $h\in\mathfrak{u}$ whenever $0\oplus hM\subseteq V$. In that case the canonical morphism $H\longrightarrow \Hom_k(M,k\oplus
M)$ in $\Mod H$ induces a monomorphism $\frac{H}{\mathfrak{a}}\rightarrowtail \Hom_k(M,\frac{k\oplus M}{V})$. But
$\Hom_k(M,\frac{k\oplus M}{V})$ is isomorphic to $D(M)^{(r)}$, for some natural number $r$, where $D=\Hom_k(?,k)$ is
the canonical duality. Then $D(M)=\tau_H^{-1}(T)$, where $T=D(S)$ is the simple projective right $H$-module. It is known that all
cyclic submodules of $D(M)$ are projective (see, e.g., \cite[Section 3.2]{Ringel}). From that one easily derives that, in
our case, $\frac{H}{\mathfrak{a}}$ is a projective right $H$-module and, hence, Êcondition 2 of Lemma \ref{new-splitting}
holds for $n=1$. However, according to Theorem \ref{clasificacion-right-split-para-buenos}, the TTF-triple in $\Mod A$ associated to $I$ is not right split
because $_HM$ is not FP-injective.

\item Let $H,C$ be rings, the first one being hereditary perfect, and $_HM_C$ be a bimodule such that $_HM$ is Êfaithful. The idempotent ideal
\[I\arr{\sim}\scriptsize{\left[\begin{array}{cc}C&0\\M&0\end{array}\right]}
\]
of
\[A\arr{\sim}\scriptsize{\left[\begin{array}{cc}C&0\\M&H\end{array}\right]}
\]
is clearly pure on the left and Ê$\op{lann}_A(I)=0$. Combining
Theorem \ref{clasificacion-right-split-para-buenos} and Theorem
\ref{clasificacion-right-split}, one gets that $I$ is right
splitting (i.e. it satisfies the equivalent conditions of Lemma
\ref{new-splitting}) if, and only if, $_HM$ is FP-injective
(equivalently, $_HM$ has a hereditary $\Pi$-projective dual). We
leave as an exercise to check it directly by using an argument
similar to 1.

Ê\item If $A$ is commutative and we denote by $\mathfrak{L}$, $\mathfrak{C}$ and $\mathfrak{R}$ the sets of left, centrally and right split
ÊTTF-triples in $\Mod A$, respectively, then Ê$\mathfrak{L}=\mathfrak{C}\subset\mathfrak{R}$ and the last inclusion may be strict.
ÊIndeed, since all idempotents in $A$ are central, the equality $\mathfrak{L}=\mathfrak{C}$ follows from Corollary \ref{clasificacion-left-split}. On the other hand, if
Ê$k$ is a field and $A=\{\lambda =(\lambda_n)\in k^\N:$ $\lambda\text{ is eventually constant}\}$, then $I=k^{(\N)}$ is an idempotent ideal of
Ê$A$ which is pure and satisfies that $\op{ann}_A(I)=0$. Moreover, one has $A/I\cong k$ and then condition $(3)$ in TheoremÊ~\ref{clasificacion-right-split} holds
Ê(see Example \ref{ejemplo de ideal right-splitting}). The associated TTF-triple is then right split but not centrally split.
\end{enumerate}
\end{remexa}

\printindex

\end{document}